\documentclass[]{svproc}

\usepackage{amsmath,amssymb}
\usepackage[ruled,linesnumbered]{algorithm2e}
\usepackage[colorlinks=true,linkcolor=blue,citecolor=blue,urlcolor=blue]{hyperref}
\usepackage{enumitem} % Provide a more flexible spacing for enumerate.
\usepackage{multirow}
\usepackage{array}
\newcolumntype{L}[1]{>{\raggedright\let\newline\\\arraybackslash\hspace{0pt}}m{#1}}
\newcolumntype{C}[1]{>{\centering\let\newline\\\arraybackslash\hspace{0pt}}m{#1}}
\newcolumntype{R}[1]{>{\raggedleft\let\newline\\\arraybackslash\hspace{0pt}}m{#1}}

\usepackage{color}
\usepackage[obeyDraft]{todonotes} % REMOVE THIS once the draft is complete.
\newcommand{\abs}[1]{\lvert#1\rvert}
\newcommand{\norm}[1]{\lVert #1 \rVert}

\newcommand{\union}{\cup}

\newcommand{\naturals}{\ensuremath{\mathbb{N}}}
\newcommand{\integers}{\ensuremath{\mathbb{Z}}}
\newcommand{\rationals}{\ensuremath{\mathbb{Q}}}
\newcommand{\reals}{\ensuremath{\mathbb{R}}}
\newcommand{\complexes}{\ensuremath{\mathbb{C}}}

\newcommand{\PSLQ}{\textsc{pslq}}
\newcommand{\APSLQ}{\textsc{apslq}}

\newcommand{\GOOD}{\textsc{good}}
\newcommand{\UNEXPECTED}{\textsc{unexpected}}
\newcommand{\BAD}{\textsc{bad}}
\newcommand{\FAIL}{\textsc{fail}}

\newcommand{\Qext}[1]{\ensuremath{\rationals[#1]}}
\newcommand{\Zext}[1]{\ensuremath{\integers[#1]}}

\newcommand{\quadint}[1]{\ensuremath{\mathcal{O}_{\rationals[#1]}}}
\newcommand{\nint}[1]{\ensuremath{\left\lceil #1 \right\rfloor}}
\newcommand{\floor}[1]{\ensuremath{\left\lfloor #1 \right\rfloor}}
\newcommand{\ceil}[1]{\ensuremath{\left\lceil #1 \right\rceil}}
\newcommand{\defeq}{\ensuremath{\stackrel{\text{\tiny{}def}}{=}}}

\DeclareMathOperator*{\argmax}{arg\,max}
\DeclareMathOperator*{\argmin}{arg\,min}
\DeclareMathOperator{\col}{col}
\DeclareMathOperator{\row}{row}

\newcommand{\orcid}[1]{\textsc{orc}i\textsc{d}: \href{https://orcid.org/#1}{\texttt{\footnotesize #1}}}

\spnewtheorem*{notation}{Notation}{\bfseries\upshape}{\itshape}

\title{Extending the PSLQ Algorithm to Algebraic Integer Relations}

\author{Matthew P. Skerritt \and Paul Vrbik}
\institute{
	Centre for Computer-assisted Research Mathematics and its Applications (\textsc{carma}),\\
	School of Mathematical and Physical Sciences,
	University of Newcastle,
	Australia,\\
	\email{matthew.skerritt@uon.edu.au},
	\orcid{0000-0003-2211-7616}
	\\
	\email{paulvrbik@gmail.com}
}

\date{\today}

\begin{document}
	\maketitle

	\begin{abstract}
		The \PSLQ{} algorithm computes integer relations for real numbers and Gaussian integer relations for complex numbers. We endeavour to extend \PSLQ{} to find integer relations consisting of algebraic integers from some quadratic extension fields (in both the real and complex cases). We outline the algorithm, discuss the required modifications for handling algebraic integers, problems that have arisen, experimental results, and challenges to further work.
	\end{abstract}
	
	\section{Introduction} % (fold)
	\label{sec:introduction}
		The Euclidean algorithm for real numbers\footnote{Euclid's Elements Book 10} is perhaps the simplest example of an \emph{integer relation algorithm}. 
		Given \( a, b, \in \reals \) the algorithm computes \( g \in \reals \) such that \( a = mg \) and \( b = ng \) for some \( m,n \in \integers\). 
		If we let \( s = n \) and \( t = -m \) then we have found the relation \( as + bt = 0 \).
		It was Ferguson and Forcade's efforts to generalise this to the case where \( a_1, \dotsc, a_n \in \reals\) in 1979 \cite{Ferguson:1979} that eventually led to the \PSLQ{} algorithm by Ferguson and Bailey in 1991 \cite{Ferguson:1991}.

		This general case is attractive. One may determine if a number \(\alpha\) is algebraic by finding an integer relation for \(\left( \alpha^0, \, \alpha^1,\, \dotsc, \alpha^n\right)\) for some \(n \in \mathbb{N}\).
		Furthermore, searching for such relations involving \(\pi\) led to the discovery of the Bailey-Borwein-Plouffe (\textsc{bbp}) formula \cite{BBP:1996}.

		A further extension of the integer relation problem is from real numbers and integers to complex numbers and Gaussian integers respectively.
		This extension was shown to be handled by the \PSLQ{} algorithm in the 1999 paper by Ferguson, Bailey and Arno \cite{Ferguson:1999} in which they analysed the algorithm and proved bounds on the number of iterations required to find a relation.
		The complex case is rarely mentioned in the literature, although we note that it is handled by \emph{Maple}'s implementation of the algorithm.
		
		The integer relation cases handled by the \PSLQ{} algorithm are covered by the following definition.
		\begin{definition}[Integer Relation]\label{def:integer_relation}
			Let \( \mathbb{F} \in \{ \reals, \complexes \} \), and let
			\[ \mathcal{O} = \begin{cases} \integers & \text{ if } \mathbb{F} = \reals \\ \Zext{\sqrt{-1}} & \text{ if } \mathbb{F} = \complexes \end{cases} \]
			For \( x \in \mathbb{F}^n \), an \emph{integer relation} of $x$ is a vector \( a \in \mathcal{O}^n \), \( a \ne 0 \), such that \(a_1 x_1 + \dots + a_n x_n = 0\).
		\end{definition}
		
		We will further generalise the integer relation problem in this paper.
		In order to talk about the algorithm in more generality we will use the following notation.
		\begin{notation}[\(\mathbb{F},\mathcal{O}\)]
			When discussing \PSLQ{} and our generalisations we will denote by \(\mathbb{F}\) the field from which the input to the algorithm is taken, and by \(\mathcal{O}\) the ring of integers from which the elements of the integer relation belong.
		\end{notation}
		
		Observe that for the linear combination property of an integer relation to be well defined, it must be the case that \(\mathcal{O} \subset \mathbb{F}\). 
		As such, we may consider the notion a nearest integer to a given element of the field.
		This is important for the \PSLQ{} algorithm.
		
		\begin{definition}[Nearest Integer]\label{def:nint}
			Let \(x \in \mathbb{F}\).  An integer \(a \in \mathcal{O}\) is a \emph{nearest integer} to \(x\) if \(\abs{x-a}\) is minimal.
			We consider a function \( \nint{\cdot} : \mathbb{F} \to \mathcal{O} \) to be a \emph{nearest integer function} if it maps each \(x \in \mathbb{F}\) to one of its nearest integers. 
			When the ring of integers needs to be specified, we will denote a nearest integer function by \( \nint{\cdot}_{\mathcal{O}} \).
		\end{definition}
		
		\subsection{Algorithm Overview} % (fold)
		\label{sub:algorithm_overview}
		We provide a high level description of the unmodified \PSLQ{} algorithm which is sufficient to understand the modifications we have made.
		For an alternative and slightly more detailed introduction the reader is referred to Straub \cite{Straub:}.
		
		We show the mathematical details of the algorithm, but omit many technical considerations needed for a practical and effective implementation.
		Details and analysis suitable for a practical implementation can be found in the literature, in particular: Borwein \cite{borwein2002computational}, and Bailey and Broadhurst \cite{Bailey:2001}.
		
		The \PSLQ{} algorithm has parameters \(\tau\), \(\gamma\), and \(\rho\) that must satisfy
		\begin{gather}
			\frac{1}{\rho} \ge \abs{x-\nint{x}} \quad \forall x \in \mathbb{F} \label{cnd:rho} \\ 
			1 < \tau \le \rho \label{cnd:tau_rho} \\
			\frac{1}{\tau^2} = \frac{1}{\gamma^2} + \frac{1}{\rho^2} \label{cnd:tau_rho_gamma}
		\end{gather}
		in order to establish runtime bounds on the algorithm \cite{Ferguson:1999}.
		
		For a given \( \mathbb{F} \), so long as \( \mathcal{O} \) is a lattice, there exists \( \rho \) such that the inequality (\ref{cnd:rho}) is sharp.
		Using this value for \( \rho \) gives the most flexibility with the other parameters.
		From (\ref{cnd:tau_rho_gamma}) we see that \(\tau \to \rho \) as \( \gamma \to \infty \) and that for fixed \( \rho \) there will be a greatest lower bound for \( \gamma \) such that \( \tau > 1 \).

		\begin{definition}[\(\gamma_1\)]\label{def:gamma_1}
			Let \(\rho \) be such that (\ref{cnd:rho}) is sharp. Then \(\gamma_1\) is the value of \(\gamma\) that satisfies \( 1 = 1/\gamma^2 + 1/\rho^2 \).
		\end{definition}

		We use the value of \( \rho \) such that (\ref{cnd:rho}) is sharp, and choose any \(\gamma > \gamma_1\). 
		So long as \( \rho \) > 1 (i.e., \(1/\rho\) < 1) then all three conditions will be satisfied. 
		
		Note that when \( \mathbb{F} = \reals \) and \( \mathcal{O} = \integers \) then the above strategy gives \( \rho = 2 \) and \( \gamma_1 = \sqrt{4/3} \). 
		This value of \( \gamma_1 \) is precisely the lower bound of \( \gamma \) given in the literature.
		
		Similarly when \( \mathbb{F} = \complexes \) and \( \mathcal{O} = \Zext{\sqrt{-1}} \) (i.e., Gaussian integers) then \( \rho = \sqrt{2} \) and \( \gamma_1 = \sqrt{2} \). This is precisely the bound on \( \gamma \) given in the literature for the complex case. 
		
		The \PSLQ{} algorithm is presented in Algorithm~\ref{alg:PSLQ}, below. In order to make sense of it, we need the following definitions.
		
		\begin{definition}[Lower Trapezoidal]\label{def:lower_trapezoidal}
			Let \( H = \left( h_{i,j} \right) \) be an \( m \times n \) matrix. If \( h_{i,j} = 0 \) whenever \( j > i \) then \( H \) is \emph{lower trapezoidal}.
		\end{definition}
		
		Note that a lower trapezoidal square matrix is exactly a lower triangular matrix.
		
		\begin{definition}[\(H_x\)]\label{def:H_x}
			Let \( x \in \mathbb{F}^n \). Then the \( n \times (n-1) \) matrix \( H_x = \left( h_{i,j} \right) \) is defined by
			\[
				h_{i,j} = 	\begin{cases}
								\,0 & \text{if } i < j \\
								\,s_{i+1}/{s_i} & \text{if } i = j \\
								\,-\overline{x_i}x_j/(s_js_{j+1}) & \text{if } i > j
							\end{cases}
				\qquad
				\text{ where }
				s_i = \sqrt{\sum_{k=i}^n x_k\overline{x_k}}
			\]
		\end{definition}
		
		Note that the complex conjugates are needed for full generality to cope with the complex case. Often the literature will present only the real case of \PSLQ{} in which case \( x_k\overline{x_k} = x_k^2 \) and is reported as such. Similarly for the conjugates in Definition~\ref{def:Q}, below.
		
		\begin{definition}[Hermite Reduction, \(D_H\)]
			Let \(A = (a_{i,j})\) be a lower trapezoidal \( m \times n \) matrix with \(a_{j,j}\ne 0 \) for all \( j \). Then the \( m \times m \) matrix \(D_A = (d_{i,j})\) where
			\[
				d_{i,j} = 	\begin{cases}
								\,0 & \text{if } i < j \\
								\,1 & \text{if } i = j \\
								\,\nint{\displaystyle\frac{-1}{a_{j,j}}\sum_{k=j}^id_{i,k}a_{k,j}} & \text{if } i > j
							\end{cases}
			\]
			is the \emph{reducing matrix} of \(A\). The matrix \(D_A\,A\) is the \emph{Hermite reduction} of \(A\).
		\end{definition}

		Observe that \(D_H\) is a lower triangular matrix containing invertible integers on its diagonal. It is therefore an invertible matrix whose inverse is also integer valued.
		
		\begin{definition}[\(Q_{[A,k]}\)]\label{def:Q}
			Let \(A = (a_{i,j})\) be an \( m \times n \) matrix with \( m > n \), and let \( 1 \le k \le n \). Let \( \beta = a_{k,k} \), \( \lambda = a_{k,k+1} \), and \( \delta = \sqrt{\beta\overline{\beta} + \lambda\overline{\lambda}} \) Then the \( n \times n \) block diagonal matrix
			\[
				Q_{[A,k]} \defeq \begin{cases} \,I_{n} & \text{if } k = n \\ \, \left(q_{i,j}\right) & \text{otherwise} \end{cases}
			\]
			where \( \left(q_{i,j}\right) \) is the block diagonal matrix with submatrix
			\[ 
				\begin{pmatrix} q_{k,k} & q_{k,k+1} \\ q_{k+1,k} & q_{k+1,k+1} \end{pmatrix}
				=
				\frac{1}{\delta}\begin{pmatrix} \,\overline{\beta} & \,-\lambda \\ \,\overline{\lambda} & \,\beta \end{pmatrix}
			\]
			and 1's for all other diagonal entries.
		\end{definition}
		
		Observe that multiplication on the right by \( Q_{[A,k]} \) changes only columns \( k \) and \(k+1\) in a way that is effectively multiplying those columns as a submatrix by the submatrix explicitly stated in the definition.
		
		When used in Algorithm~\ref{alg:PSLQ} (line~\ref{PSLQ:corner}) \( Q_{[H^\prime,r]} \) is an orthogonal matrix. The swapping of rows that occurs in the prior steps will usually cause \( H^\prime \) to cease to be lower trapezoidal. The post-multiplication with \(  Q_{[H^\prime,r]} \) ensures that \( H^\prime \) is once again lower trapezoidal \cite{borwein2002computational,Ferguson:1999}. The only case where the row swap does not remove the lower-trapezoidal property of \( H^\prime \) is when \( r = n-1 \) in which case \( Q_{[H^\prime,r]} \) is the identity matrix and so \( H^\prime \) is unaffected.
		
		Finally, we use the following notation to refer to rows and columns of matrices, when needed.
		
		\begin{notation}[\(\col_k, \row_k\)]\label{def:col_row}
			For a matrix \( M \) we denote by \(\col_k(M)\) the \(k\)th column of \(M\) and by \(\row_k(M)\) the \(k\)th row of M.
		\end{notation}

		\begin{algorithm}
			\caption{\PSLQ}
			\label{alg:PSLQ}
			\DontPrintSemicolon
			\SetKwInOut{Input}{input}
			\SetKwInOut{Output}{output}
			\Input{\(x \in \mathbb{F}^n, \gamma > \gamma_1\)}
			\Output{\( a \in \mathcal{O}^n \)}
			\SetKwRepeat{Do}{do}{while}
			\SetKwRepeat{Repeat}{repeat}{until}
			\SetKw{KwAnd}{and}
			\SetCommentSty{commentstyle}
			\BlankLine
			\tcc{\rule[0.5ex]{0.3\textwidth}{1px} \hfill Initialisation \hfill \rule[0.5ex]{0.3\textwidth}{1px}}
			\( H^\prime \leftarrow H_{x / \norm{x}} \quad A \leftarrow I_n \)\;
			\BlankLine
			\tcc{\rule[0.5ex]{0.3\textwidth}{1px} \hfill Main Calculation \hfill \rule[0.5ex]{0.3\textwidth}{1px}}
			\Repeat{\(r = n-1\) \KwAnd \( H^\prime_{n-1,n-1} = 0\)} {
				\(H^\prime \leftarrow D_{H^\prime}\,H^\prime\)\tcc*{Hermite reduce \(H^\prime\)}
				\(A \leftarrow D_{H^\prime}\,A\)\tcc*{Update \(A\)}
				\(r \gets \argmax_{1\le r \le n-1}(\gamma^r\abs{H^\prime_{r,r}})\)\tcc*{Find \(r\) such that \(\gamma^r\abs{H^\prime_{r,r}}\) is maximal}
				\( \row_r(H^\prime) \leftrightarrow \row_{r+1}(H^\prime) \)\tcc*{Exchange rows \(r\) and \(r+1\) in \(H^\prime\)}
				\( \row_r(A) \leftrightarrow \row_{r+1}(A) \)\tcc*{Exchange rows \(r\) and \(r+1\) in \(A\)}
				\( H^\prime \leftarrow H^\prime\,Q_{[H^\prime,r]} \)\nllabel{PSLQ:corner}  \tcc*{Make sure \(H^\prime\) is lower trapezoidal}
			}
			\BlankLine
			\KwRet{\(\col_n(A^{-1})\)}\;
		\end{algorithm}
		
		After each iteration the value \( 1 / \max{\abs{H^\prime_{r,r}}} \) is a lower bound for \emph{any} integer relation of \( x \). 
		Furthermore if \( a \) is the integer relation found by the algorithm, then \( \norm{a} \le \gamma^{n-2} M\) where \(M\) is the norm of the smallest possible integer relation \cite[Theorem 3]{Ferguson:1999}.
		
		Note that the algorithm as presented above does not terminate if there is no integer relation for the input \( x \). 
		This can be remedied either by specifying termination after a maximum number of iterations are performed, or after the lower bound for the norm of an integer relation exceeds some value.
		
		The algorithm is exact if the individual steps can be performed exactly. 
		That is to say, if we could compute with all real numbers exactly then the algorithm would always calculate an integer relation if there is one to be found.
		Furthermore, it will find an integer relation in a polynomially bounded number of iterations \cite{borwein2002computational,Ferguson:1999}.
		In practice, however, an implementation of the \PSLQ{} algorithm must use floating point arithmetic and so numerical error may prevent the detection of a valid integer relation.
		Nonetheless \PSLQ{} has shown remarkable numerically stability.
		
		Finally, we reiterate that the algorithm as presented here lacks the details needed for practical numeric application. There are many optimisations that can, and should, be implemented in order for an implementation to be effective. The interested reader should consult the literature \cite{Bailey:2001,borwein2002computational}.
		% subsection algorithm_overview (end)

		\subsection{Algebraic Number Theory} % (fold)
		\label{sub:algebraic_number_theory}
			We introduce only enough algebraic number theory as is needed. 
			The reader is referred to the literature for a more thorough study \cite[e.g.]{stewart2001algebraic}.

			\begin{definition}[Algebraic Number]
				A number \(\alpha \in \complexes\) is an \emph{algebraic number} (or simply \emph{algebraic}) if it is a zero of a polynomial with rational coefficients.
			\end{definition}
	
			\begin{definition}[Algebraic Integer]
				A number \(\alpha \in \complexes\) is an \emph{algebraic integer} if it is a zero of a monic polynomial with integer coefficients.
				The ring of all algebraic integers is denoted by \( \mathcal{A} \).
			\end{definition}
	
			\begin{definition}[Algebraic Extension]
				A field, \(\mathbb{K} \supset \rationals\), is an \emph{algebraic extension field} (or simply an \emph{algebraic extension}) if \(k\) is algebraic for all \(k \in \mathbb{K}\). 
			\end{definition}
	
			We may now talk of the algebraic integers of a particular algebraic extension field.
	
			\begin{definition}\label{def:algebraic_integer_ring}
				Let \(\mathbb{K}\) be an algebraic extension field.
			    The \emph{ring of integers of \(\mathbb{K}\)}, denoted \( \mathcal{O}_\mathbb{K} \), is the intersection \( \mathbb{K} \cap \mathcal{A} \) of the extension field with the ring of all algebraic integers.
			\end{definition}
	
			For the purposes of this paper we consider only simple quadratic extension fields.
			That is, fields of the form
			\( \Qext{\sqrt{D}} := \{ q_1 + q_2\sqrt{D}\;\vert\; q_1, q_2 \in \rationals \} \).
			Without loss of generality we may assume \(D \in \integers\) is square free.
			The ring of integers of such fields are known \cite{stewart2001algebraic} to be
			\( \quadint{\sqrt{D}} = \Zext{\omega} = \left\{ \alpha + \beta\,\omega \;\vert\; \alpha, \beta \in \integers \right\} \)
			where
			\begin{equation}
				\omega =	\begin{cases}
								\sqrt{D} & \text{if } D \equiv 2,3 \text{ (mod } 4 \text{)} \\
								(1 + \sqrt{D})/2 & \text{if } D \equiv 1 \text{ (mod } 4 \text{)}
							\end{cases}
			\end{equation}
		% subsection algebraic_number_theory (end)
	% section introduction (end)

	\section{Extension to Algebraic Integers} % (fold)
	\label{sec:algebraic_pslq}
		In order to extend \PSLQ{} to allow for algebraic integers, we first establish the relationship between algebraic integers, algebraic extension fields, and integer relations. We want to generalise, and thus wish to encapsulate the cases already handled by the existing theory.
		
		A naïve strategy would be to replace \( \mathbb{F} \) in Definition~\ref{def:integer_relation} with an arbitrary extension field, and to replace \( \mathcal{O} \) with the ring of integers of that extension field.
		However, observe that the integers (\integers) are not the ring of integers of the field of real numbers.
		Similarly, the Gaussian integers (\Zext{\sqrt{-1}}) are not the ring of integers of the field of complex numbers.
		So this strategy will not capture the pre-existing cases.
		
		We instead generalise by introducing an intermediate extension field, according to the following definition.
	
		\begin{definition}[Algebraic Integer Relation]
			\label{def:algebraic_integer_relation}
			Let \(x \in \mathbb{F}^n\) and \(\mathbb{K} \subseteq \mathbb{F}\) be an algebraic extension field.
			An \emph{algebraic integer relation} of \(x\) is a vector  \(a \in \left(\mathcal{O}_\mathbb{K}\right)^n\), \( a \ne 0 \), such that \(a_1 x_1 + \dots + a_n x_n = 0\).
		\end{definition}
	
		Observe that algebraic integer relations are indeed a generalisation of integer relations.
		When \(\mathbb{F} = \reals\) and \(\mathbb{K}=\rationals\) (thinking of \rationals{} as a trivial extension field) then an algebraic integer relation is also an integer relation satisfying Definition~\ref{def:integer_relation}. The same is true for the complex case when \(\mathbb{F}=\complexes\) and \(\mathbb{K}=\Qext{\sqrt{-1}}\)).
		
		Since we have stated above that we are only concerning ourselves with simple quadratic extension fields, we correspondingly restrict our attention to algebraic integer relations where \( \mathbb{K} = \Qext{\sqrt{D}} \) is a quadratic extension field, and \( \mathbb{F} \) is the Archimedean norm closure of \( \mathbb{K} \) (i.e., \reals{} if \( D \ge 0 \) and \complexes{} if \( D < 0 \)).
			
		\subsection{Reduction} % (fold)
		\label{sub:reduction}
			One approach to computing algebraic integer relations is to reduce the problem to an integer relation problem.
			We may then solve the problem with an existing integer relation finding algorithm, such as \PSLQ.
			
			Observe that for \(\alpha + \beta\,\omega \in \quadint{\sqrt{D}}\) we have \(\left(\alpha + \beta\,\omega\right)\,x = \alpha\,x + \beta\left(x\,\omega\right)\).
			This suggests a method of reduction.
	
			Given an algebraic extension field \(\Qext{\sqrt{D}} \subset \mathbb{F}\), and input
			\(\left(x_1,\dotsc,x_n\right) \in \mathbb{F}^n \)
			we compute
			\( \left( x_1, x_1\omega, \dotsc, x_n, x_n\omega \right) \)
			which we give as input to \PSLQ{} producing an integer relation
			\( \left(a^\prime_1, \dotsc, a^\prime_{2n} \right) \)
			from which we attempt to reconstruct an algebraic integer relation
			\( \left( a_1, \dotsc, a_n \right) \) where \(a_k = a^\prime_{2k-1}+a^\prime_{2k}\omega \).
	
			When \(\mathbb{F} = \reals\) it is straightforward to see that each \(a_k \in \quadint{\sqrt{D}}\), and so the reconstructed relation is, indeed, an algebraic integer relation.
			
			However, when \(\mathbb{F}=\complexes\) the \(a^\prime_k\) are Gaussian integers \(\alpha_k + \beta_k\,i\) where \(\alpha_k,\beta_k\in\integers\).
			Then
			\begin{align*} 
				a_k &= \left(\alpha_{2k-1}+\beta_{2k-1}\,i\right) + \left(\alpha_{2k}+\beta_{2k}\,i\right)\,\omega
				    = \left(\alpha_{2k-1}+\alpha_{2k}\,\omega\right) + \left(\beta_{2k-1} + \beta_{2k}\,\omega \right)i
			\end{align*}
			which will not always be an algebraic integer in \quadint{\sqrt{D}}.
			% In general the recovered relation will consist of integers in \Qext{\sqrt{\abs{D}}+i}.
		
			Ideally, we want the \(a^\prime_k\) to only ever be integer valued.
			In some cases it may be possible to transform \(\left(a^\prime_1,\dotsc,a^\prime_n\right)\) into an equivalent (for the purposes of algebraic integer relation detection) integer valued vector, such as dividing by a common Gaussian integer divisor.
			We have not yet found a reliable way to detect such cases in general.
		% subsection reduction (end)
	
		\subsection{Algebraic PSLQ} % (fold) % NOTE: Cannot use \PSLQ in the title because of font restrictions. Smallcaps won't work in there.
			\label{sub:algebraic_pslq}
			An alternative approach to computing algebraic integer relations is to modify the \PSLQ{} algorithm to compute them directly.
			We call this modified algorithm \emph{Algebraic \PSLQ}, or \APSLQ. 
	
			We observe that the reducing matrix is the source of integers in the algorithm.
			The reducing matrix, in turn, relies on the nearest integer function.
			The theorems bounding the number of iterations needed to find an integer relation rely only on the \( \tau, \rho \), and \( \gamma \) parameters, the latter of which is arbitrarily chosen and the others of which are determined by the properties of the integer lattice.
			
			In order to utilise as much of the existing theory as possible we replace the nearest integer function in the computation of the reducing matrix with a nearest algebraic integer function. 
			Additionally, we require the specification of the intermediate quadratic extension field as input to the algorithm.
			The algorithm remains otherwise unmodified.
	
			This immediately causes a problem.
			In the case of a real quadratic extension field (when \(D > 0\)) the algebraic integers are dense in \reals.
			This leaves us without a well defined nearest integer, and hence no integer lattice.
			We put this case away pending further algorithmic modifications and restrict our attention to complex quadratic extension fields \( D < 0 \).
	
			In order to calculate the nearest integer for an arbitrary \( z \in \complexes \), we first re-write \( z = \alpha + \beta\,\omega\) and use \( \alpha \) and \( \beta \) to compute \nint{z}. 
			There are two cases.
			
			When \(  D \equiv 2,3 \text{ (mod } 4 \text{)} \), then \( \alpha = \Re(z) \) and \( \beta = \Im(z) / \sqrt{\abs{D}} \). 
			We have
			\[ \nint{z} = \nint{\alpha}_\integers + \nint{\beta}_\integers\omega \]
			
			When \( D \equiv 1 \text{ (mod } 4 \text{)} \), then \( \beta = 2\Im(z)/\sqrt{\abs{D}} \) and  \( \alpha = \Re(z)-\beta/2 \). 
			We have two candidates for the nearest integer and choose the one which is closest to \( z \).
			\[
				\nint{z} = \nint{\alpha}_\integers + \floor{\beta}\omega
				\quad\text{ or }\quad
				\nint{z} = \nint{\alpha+1/2}_\integers + \ceil{\beta}\omega
			\]
	
			We bound \(\abs{ z - \nint{z} } \le \epsilon\) for all \(z \in \complexes\) using the geometric properties of the lattices.
			\[ 
				\epsilon =	\begin{cases}
								\frac{1}{2}\sqrt{\abs{D}+1} & \text{if } D \equiv 2,3 \text{ (mod } 4 \text{)} \\
								\frac{1}{4}{\frac {\abs{D}+1}{\sqrt {\abs{D}}}} & \text{if } D \equiv 1 \text{ (mod } 4 \text{)} 
							\end{cases}
			\]
			And so we can compute the corresponding value of \(\rho\)
			\[ 
				\rho =	\begin{cases}
								\frac{2}{\sqrt{\abs{D}+1}} & \text{if } D \equiv 2,3 \text{ (mod } 4 \text{)} \\
								\frac{4\sqrt{\abs{D}}}{\abs{D}+1} & \text{if } D \equiv 1 \text{ (mod } 4 \text{)}
							\end{cases}
			\]

			However, as \(\abs{D}\) increases, the value of \(\rho\) decreases, and eventually \(\rho < 1\) making it impossible to satisfy condition (\ref{cnd:tau_rho}), and causing \( \gamma_1 \) to become complex.
			This leaves us with \(D=-2\), \(D=-3\), \(D=-7\), and \(D=-11\) as the only values of \( D \) for which the existing theory holds.
			
			We will see that even when the conditions do not hold the algorithm can still be effective (see Section~\ref{sub:complex_quadratic_extension_fields}, Table~\ref{tbl:complex_quadratic-no_gamma_1-extensions}).
			
			In this paper we examine the efficacy \APSLQ{} and the reduction method.
			We leave, for now, the question of additional modifications which may handle the problems described above.
		% subsection algebraic_pslq (end)
	% section algebraic_pslq (end)

	\section{Experimental Results} % (fold)
	\label{sec:experimental_results}
		We tested the efficacy of the above two methods experimentally. 
		To do so we used \emph{Maple}'s native \PSLQ{} implementation for reduction, and our own implementation of \APSLQ{} (written in \emph{Maple}).
		The code and results are available through GitHub\footnote{Repository: \url{https://github.com/matt-sk/Algebraic-PSLQ} Tag: JBCC}.
		
		Our implementation of \APSLQ{} is described in Algorithm~\ref{alg:APSLQ}. 
		Recall that for \APSLQ{} the matrices \(D_{H^\prime}\) are constructed using an algebraic nearest integer function.
				
		\begin{algorithm}
			\caption{\APSLQ}
			\label{alg:APSLQ}
			\DontPrintSemicolon
			\SetKwInOut{Input}{input}
			\SetKwInOut{Output}{output}
			\Input{\(x \in \mathbb{F}^n, D \in \integers, \gamma \ge 0, \epsilon > 0, max_i > 0 \)}
			\Output{A vector in \( \mathcal{O}_\mathbb{K}^n \) (where \(\mathbb{K}=\quadint{\sqrt{D}}\)) or FAIL}
			\SetKwRepeat{Do}{do}{while}
			\SetKwRepeat{Repeat}{repeat}{until}
			\SetCommentSty{commentstyle}
			\SetKw{KwBy}{by}
			\SetKw{KwOr}{or}
			\BlankLine
			\tcc{\rule[0.5ex]{0.3\textwidth}{1px} \hfill Initialisation \hfill \rule[0.5ex]{0.3\textwidth}{1px}}
			\( y \leftarrow x/\norm{X} \) \tcc*{Normalise input vector}
			\( H^\prime \gets D_{H_y}\,H_y \quad B \gets D_{H_y}^{-1} \quad y \gets y\,D_{H_y}^{-1} \)\tcc*{Initial Hermite reduction}
			\( i \gets 0 \) \tcc*{Loop counter}
			\BlankLine
			\tcc{\rule[0.5ex]{0.3\textwidth}{1px} \hfill Main Calculation \hfill \rule[0.5ex]{0.3\textwidth}{1px}}
			\Repeat{\(y_k / \norm{\col_k(B)} < \epsilon \) \KwOr \( i > max_i \)}{
				\(r \gets \argmax_{1\le r \le n-1}\left(\gamma^r\abs{H^\prime_{r,r}}\right)\) \tcc*{Find \(r\) s.t. \(\gamma^r\abs{H^\prime_{r,r}}\) is maximal}
				\( \row_r(H^\prime) \leftrightarrow \row_{r+1}(H^\prime) \) \tcc*{Swap rows \(r\) and \(r+1\) in \(H^\prime\)}
				\( \col_r(B) \leftrightarrow \col_{r+1}(B) \) \tcc*{Swap columns \(r\) and \(r+1\) in \(B\)}
				\( y_r \leftrightarrow y_{r+1} \) \tcc*{Swap elements \( r \) and \( r+1 \) in \( y \)}
				\( H^\prime \leftarrow H^\prime\,Q_{[H^\prime,r]} \) \tcc*{Make sure \(H^\prime\) is lower trapezoidal}
				\( H^\prime \gets D_{H^\prime}\,H^\prime \) \tcc*{Hermite reduce \(H^\prime\)}
				\( B \gets B\,D_{H^\prime}^{-1} \quad y \gets y\,D_{H^\prime}^{-1} \)\tcc*{Update \(B\) and \(y\)}
				\( k \gets \argmin_{1 \le k \le n}(\abs{y_k}) \) \tcc*{Find \(k\) s.t. \(\abs{y_k}\) is minimal}
				\( i \gets (i+1) \)\tcc*{Increment loop counter}
			}
			\BlankLine
			\leIf{\(y_k / \norm{\col_k(B)} < \epsilon \)} {\Return{\(\col_k(B)\)}} {\Return{FAIL}}
		\end{algorithm}
		
		The particulars are a little different to the algorithm presented in Section~\ref{sub:algorithm_overview} (Algorithm~\ref{alg:PSLQ}). 
		It is effectively the algorithm as described by Borwein \cite[fig. B.5]{borwein2002computational}, although we note that our implementation correctly handles the complex case as described above, whereas the algorithm given by Borwein is specialised to the real case.
		
		To understand the differences, first note that the matrix \(B\) is simply the matrix \(A^{-1}\) from Algorithm~\ref{alg:PSLQ}. 
		Each column of \(B\) is considered a possible integer relation of \(x\), and the vector \(y\) is kept updated so that \(y = (x/\norm{x})\,B \).
		As such, if \(y_k = 0\) for some \(k\), then \(\col_k(B)\) must be an integer relation for \(x/\norm{x}\) and thus also for \(x\). 
		We terminate if we find such a relation, or if we exceed a specified number of iterations.
		This relation, \(a\) say, will not necessarily have the properly \( \norm{a} \le \gamma^{n-2} M\) that is guaranteed for a relation given by Algorithm~\ref{alg:PSLQ}, however.
		
		Note that because we are performing numeric (floating point) computations we are unlikely to exactly compute a \(0\) element in \(y\).
		To detect termination, therefore, we consider only the smallest \(\abs{y_k}\) as the best candidate for a linear combination, and look to see if it is sufficiently close to \(0\) (i.e., less than some threshold \(\epsilon\)). 
		We scale the value of \( \abs{y_k} \) by \( \norm{\col_k(B)} \) in order to avoid missing a possible relation if the norm the column of \(B\) is particularly large. 
		For more details, the reader should refer to Borwein \cite[appendix 1]{borwein2002computational}.
		
		\subsection{Methodology} % (fold)
		\label{sub:methodology}
		
		We created collections of instances of algebraic integer problems.
		Each collection, referred to as a \emph{test set}, consisted of 1000 algebraic integer relation problems. 
		
		For each test set we chose a quadratic extension field \( \mathbb{K} \), a set of constants from which we created each of the individual problems within the set, and a size for the coefficients of any algebraic integers used as part of the individual problem creation.
		
		We will speak of the choice of extension field in more detail when we describe the results, below.
		
		Two sets of constants were used: one containing real constants, the other complex. The real set was
		\begin{gather*} 
			\left\{ \pi^k : k \in \naturals, k \le 9 \right\} \union  \left\{ e^k : k \in \naturals, k \le 9 \right\} \union  \left\{ \gamma^k : k \in \naturals, k \le 9 \right\} \\
			\union  \left\{ \sin{k} : k \in \naturals, k \le 9 \right\} \union  \left\{ \log{2}, \log{3}, \log{5}, \log{7} \right\} 
		\end{gather*}
		The complex set was generated by randomly choosing an integer modulus between \( 1 \) and \( 9 \) for each integer argument from \( -9 \) to \( 9 \).
		\begin{gather*} 
			\left\{ 5\,{e^{-9\,i}},4\,{e^{-8\,i}},9\,{e^{-7\,i}},5\,{e^{-6\,i}},2\,{e^{-5\,i}},9\,{e^{-4\,i}},8\,{e^{-3\,i}},3\,{e^{-2\,i}},2\,{e^{-i}},\right.\\
			\left.4,4\,{e^{i}},5\,{e^{2\,i}},2\,{e^{3\,i}},7\,{e^{4\,i}},6\,{e^{5\,i}},3\,{e^{6\,i}},3\,{e^{7\,i}},5\,{e^{8\,i}},5\,{e^{9\,i}} 
			\right\} 
		\end{gather*}
		Each constant set was used in multiple test sets. 
		
		The size of the coefficients of the algebraic integers fell into two cases: \emph{small} (coefficients in the range \([-9,9]\) thus having exactly 1 decimal digit) and \emph{large} (coefficients in the range \([-999999,999999]\) thus having up to 6 decimal digits).
		
		Once the above choices were made for a particular test set, the problems within that set were randomly generated as follows:
		\begin{enumerate}
			\item Randomly choose an integer \( 2 \le k \le 10 \).
			\item Randomly choose \(k\) constants, \(C_1, \dotsc, C_k \), from the set of constants for the test set.
			\item For each \(C_i\), randomly choose integers \( \alpha_i \) and \(\beta_i\) within the specified size. Let \(z_i = \alpha_i + \beta_i\,\omega\).
			\item Let \(C_0=\sum_{i=1}^k z_i C_i\).
		\end{enumerate}
		The problem instance was the input vector \( x = \left(C_0, C_1, \dotsc, C_k\right) \) which, by construction, had algebraic integer relation  \( \mathfrak{a} = \left( -1, z_1, \dotsc, z_k \right) \).

		For each test set, we attempted to solve the problems within it using \PSLQ{}, reduction, and/or \APSLQ{} as appropriate.
		Our aim was to see if the algorithm could recover the known algebraic integer relation from the input vector.
		Any algebraic integer multiple of the known relation was considered to be an equivalent relation for this purpose.
		
		Test sets that used small coefficients for algebraic integers were tested using \(75\) decimal digits of floating point precision. 
		Test sets that used large coefficients were tested using \(175\) decimal digits of floating point precision.
	
		The result of a computation on an individual test instance was classified as outlined in Table~\ref{tbl:results_classifications}.
		We simply counted the number of occurrences of each result.
	
		\begin{table}
			\caption{Result Classifications}
			\begin{center}
				\begin{tabular}{ll}
					\hline
					\GOOD & The generated algebraic integer relation was recovered. \\
					\UNEXPECTED & A different, correct algebraic integer relation was found. \\
					\BAD & An incorrect algebraic integer relation was found.\\
					\FAIL & The algorithm produced no result. \\ 
					\hline
				\end{tabular}
			\end{center}
			\label{tbl:results_classifications}     
		\end{table}
		
		No \UNEXPECTED{} results were found during our testing. 
		This classification was originally introduced in the testing of an early implementation as a result of an oversight in which \(\log{2}, \log{3}\) and \(\log{6}\) were together in some problems. 
		This oversight has since been corrected, yet it remains possible (although unlikely) that other unexpected relations may still be computed, so we keep the classification as a possibility.
	
		To assess each result classification, we first note that a \FAIL{} condition is immediate if no result is produced (usually because the maximum number of iterations was exceeded).
		Assuming this is not the case, let \( a = \left( a_1, \dotsc, a_n \right) \) be the computed algebraic integer relation. 
		Let \( \mathfrak{a} = \left( -1, z_1, \dotsc, z_k \right) \) be the known relation from above. 
		Recall that we are considering any algebraic multiple of \(\mathfrak{a}\) to be correct and observe that if \( a = \lambda\mathfrak{a} \) then it must be the case that \( \lambda = -a_1 \).
		We therefore look to see if \( (-a_1)\mathfrak{a} = a \), and if so we diagnose a \GOOD{} result.
		If that is not the case, we then test the computed algebraic integer relation to \( 1000 \) decimal digits of precision, and if the result is within \(10^{-998}\) of \( 0 \) we diagnose an \UNEXPECTED{} result.
		If none of the above apply, then we diagnose a \BAD{} result.
		
		Observe that the problem with the reconstructed relation for the reduction method in the complex case as described in Section~\ref{sub:reduction} is not addressed at all by this diagnosis method.
		It is entirely possible that \( (-a_1)\mathfrak{a} = a \) even if \( a_1 \) is not a valid algebraic integer for the extension field in question.
		We describe how we accounted for this below (see Section~\ref{sub:complex_quadratic_extension_fields}).
		
		When testing sets appropriate for \APSLQ{} we performed each computation multiple times with different values of \( \gamma \) and different thresholds  for detecting integer relations in \( A^{-1} \) (as described above).
		Specifically, we used \(\gamma = \gamma_1\), \(\gamma = 2.0\), and \(\gamma =3.0\).
		Note that although the strict conditions from Section~\ref{sub:algorithm_overview} require \( \gamma > \gamma_1 \) the choice of \( \gamma = \gamma_1 \) seems to be common in practice, and the results below do not seem to suffer.
		
		The thresholds used were  \( 10^{-(d-1)}, 10^{-(d-4)}, \) and \( 10^{-(d-\log_{10}{n})} \) where \(d\) is the floating point precision in decimal digits, and \(n\) is the number of elements in the input vector.
		Note that the latter of these, copied from \emph{Maple}'s implementation, varies slightly with the number of elements of the input vector.
		These different thresholds made almost no difference whatsoever. 
		For the cases where there is no \( \gamma_1 \) (see Table~\ref{tbl:complex_quadratic-no_gamma_1-extensions}) the latter threshold sometimes had one fewer \GOOD{} and one more \FAIL{} result when compared to the other thresholds.
		We do not consider this significant, and report the results for the first threshold (\(10^{-(d-1)}\)) only.
		
		The test sets fell into three broad categories, described separately in the subsections that follow.
		% subsection subsection_name (end)
				
		\subsection{Real and Complex PSLQ} % (fold)
		\label{sub:real_and_complex_pslq}
		\todo[inline]{Section needs a better name}
		We tested our implementation of \APSLQ{} against \emph{Maple}'s \PSLQ{} implementation for the cases that \PSLQ{} was already known to work for. 
		That is for the trivial case \( \mathbb{K} = \rationals \) and the case \( \mathbb{K} = \Qext{\sqrt{-1}} \). 
		This testing acted as a “sanity check” that our implementation was correct in the known cases.
		
		\begin{table}
			\centering
			\caption{Direct comparison of \PSLQ{} and \APSLQ{}}
			\begin{tabular}{ccccccccccc}
				\hline
				\multirow{3}{*}{Field} & \hspace{2ex} & \multicolumn{4}{c}{Small Coefficients} & \hspace{2ex} & \multicolumn{4}{c}{Large Coefficients} \\
				& & \multirow{2}{*}{\PSLQ} & \multicolumn{3}{c}{Algebraic \PSLQ} & & \multirow{2}{*}{\PSLQ} & \multicolumn{3}{c}{Algebraic \PSLQ} \\
				& & & \( \gamma=\gamma_1 \) & \( \gamma = 2.0 \) & \( \gamma = 3.0 \) & & & \( \gamma=\gamma_1 \) & \( \gamma = 2.0 \) & \( \gamma = 3.0 \)\\
				\hline
				\\
				\multicolumn{11}{l}{Real \(C_i\)} \\
				\hline
				\rule{0pt}{2.5ex}\rationals & & \textsc{1000g} & \textsc{1000g} & \textsc{1000g} & \textsc{1000g} & & \textsc{1000g} & \textsc{1000g} & \textsc{1000g} & \textsc{1000g} \\[0.25ex]
				\Qext{\sqrt{-1}} & & \textsc{1000g} & \textsc{1000g} & \textsc{1000g} & \textsc{1000g}& & \textsc{1000g} & \textsc{1000g} & \textsc{1000g} & \textsc{1000g} \\[0.5ex]
				\hline
				\\
				\multicolumn{11}{l}{Complex \(C_i\)} \\
				\hline
				\rule{0pt}{2.5ex}\Qext{\sqrt{-1}} & & \textsc{1000g} & \textsc{1000g} & \textsc{1000g} & \textsc{1000g}& & \textsc{1000g} & \textsc{1000g} & \textsc{1000g} & \textsc{1000g} \\[0.5ex]
				\hline
			\end{tabular}
			\label{tbl:sanity_check}
		\end{table}

		The results are tabulated in Table~\ref{tbl:sanity_check}. Note that it is impossible to create test sets that use complex \( C_i \) and \( \mathbb{K} = \rationals \), so we were only able to test a single field with complex constants.
        % subsection real_and_complex_pslq (end)
		
		\subsection{Real Quadratic Extension Fields} % (fold)
		\label{sub:real_quadratic_extension_fields}
		
		For the real quadratic algebraic integer relations we tested the following real quadratic extension fields:
		\[ \mathbb{K} = \Qext{\sqrt{D}} \quad\text{for}\quad D \in \{2,3,5,6,7,10,11\} \]
		Recall that \APSLQ{} is not appropriate for these extension fields, so only the reduction method was tested.
		
		\begin{table}
			\centering
			\caption{Real quadratic fields, Real \( C_i \)}
			\begin{tabular}{ccccccc}
				\hline
				\multirow{2}{*}{Field} & \hspace{2ex} & \multicolumn{2}{c}{Small Coefficients} & \hspace{2ex} & \multicolumn{2}{c}{Large Coefficients} \\
				& & Reduction & \APSLQ & & Reduction & \APSLQ \\
				\hline
	            \rule{0pt}{2.5ex}\Qext{\sqrt{2}} & & \textsc{1000g} & \textsc{n/a} & & \textsc{1000g} & \textsc{n/a} \\[0.25ex]
	            \Qext{\sqrt{3}} & & \textsc{1000g} & \textsc{n/a} & & \textsc{1000g} & \textsc{n/a} \\[0.25ex]
	            \Qext{\sqrt{5}} & & \textsc{1000g} & \textsc{n/a} & & \textsc{1000g} & \textsc{n/a} \\[0.25ex]
	            \Qext{\sqrt{6}} & & \textsc{1000g} & \textsc{n/a} & & \textsc{1000g} & \textsc{n/a} \\[0.25ex]
	            \Qext{\sqrt{7}} & & \textsc{1000g} & \textsc{n/a} & & \textsc{1000g} & \textsc{n/a} \\[0.25ex]
	            \Qext{\sqrt{10}} & & \textsc{1000g} & \textsc{n/a} & & \textsc{1000g} & \textsc{n/a} \\[0.25ex]
	            \Qext{\sqrt{11}} & & \textsc{1000g} & \textsc{n/a} & & \textsc{1000g} & \textsc{n/a} \\[0.5ex]
				\hline
			\end{tabular}
			\label{tbl:real_quadratic}
		\end{table}	
		
		The results are tabulated in Table~\ref{tbl:real_quadratic}. 
		We note that since we are testing real quadratic extension fields we are in the case where \( \mathbb{F} = \reals \) and so, as stated in Section~\ref{sub:reduction}, we definitely have found algebraic integer relations. 
		Contrast this to the complex quadratic extension field testing, below.
		% subsection real_quadratic_extension_fields (end)
		
		\subsection{Complex Quadratic Extension Fields} % (fold)
		\label{sub:complex_quadratic_extension_fields}
		
		For the real quadratic algebraic integer relations we were able to test both the reduction method, and \APSLQ{}.
		We tested the following complex quadratic extension fields:
		\[ \mathbb{K} = \Qext{\sqrt{D}} \quad\text{for}\quad D \in \{-2,-3,-5,-6,-7,-10,-11\} \]
		As we have tested both reduction and \APSLQ{} for these fields, we may compare the relative efficacy of the two methods. 
		
		We accounted for the reduction problem described in Section~\ref{sub:reduction} by checking to see if the entries in the recovered relation consisted only of valid algebraic integers from the appropriate field. 
		This check was performed after the usual diagnosis, so that we could compare these \FAIL{} results with the originally diagnosed result.
		If any entries were not appropriate algebraic integers then we changed the diagnosed result to a \FAIL{} and also recorded the old result.
		We note that all such \FAIL{} results reported for our reduction tests were initially \GOOD{} results.
				
		The cases where \( D \in \{-2,-3,-7,-11\} \) are cases where \( \gamma_1 \) exists and so the three conditions (\ref{cnd:rho}), (\ref{cnd:tau_rho}), and (\ref{cnd:tau_rho_gamma}) from Section~\ref{sub:algorithm_overview} are satisfied.
		These results are summarised in Table~\ref{tbl:complex_quadratic-gamma_1-extensions}.
		Both reduction and \APSLQ{} perform superbly for these cases.

	\begin{table}
		\centering
		\caption{Complex quadratic fields with \( \gamma_1 \)}
		\begin{tabular}{ccccccccccc}
			\hline
			\multirow{3}{*}{Field} & \hspace{2ex} & \multicolumn{4}{c}{Small Coefficients} & \hspace{2ex} & \multicolumn{4}{c}{Large Coefficients} \\
			& & \multirow{2}{*}{Reduction} & \multicolumn{3}{c}{Algebraic \PSLQ} & & \multirow{2}{*}{Reduction} & \multicolumn{3}{c}{Algebraic \PSLQ} \\
			& & & \( \gamma=\gamma_1 \) & \( \gamma = 2.0 \) & \( \gamma = 3.0 \) & & & \( \gamma=\gamma_1 \) & \( \gamma = 2.0 \) & \( \gamma = 3.0 \)\\
			\hline
			\\
			\multicolumn{11}{l}{Real \(C_i\)} \\
			\hline
			\rule{0pt}{2.5ex}\Qext{\sqrt{-2}} & & \textsc{912g88f} & \textsc{1000g} & \textsc{1000g} & \textsc{1000g} & & \textsc{952g48f} & \textsc{1000g} & \textsc{1000g} & \textsc{1000g} \\[0.25ex]
			\Qext{\sqrt{-3}} & & \textsc{919g81f} & \textsc{1000g} & \textsc{1000g} & \textsc{1000g} & & \textsc{923g77f} & \textsc{1000g} & \textsc{1000g} & \textsc{1000g} \\[0.5ex]
			\Qext{\sqrt{-7}} & & \textsc{956g44f} & \textsc{1000g} & \textsc{1000g} & \textsc{1000g} & & \textsc{949g51f} & \textsc{1000g} & \textsc{1000g} & \textsc{1000g} \\[0.5ex]
			\Qext{\sqrt{-11}} & & \textsc{975g25f} & \textsc{1000g} & \textsc{1000g} & \textsc{1000g} & & \textsc{981g19f} & \textsc{1000g} & \textsc{1000g} & \textsc{1000g} \\[0.5ex]
			\hline
			\\
			\multicolumn{11}{l}{Complex \(C_i\)} \\
			\hline
			\rule{0pt}{2.5ex}\Qext{\sqrt{-2}} & & \textsc{911g1b88f} & \textsc{1000g} & \textsc{1000g} & \textsc{1000g} & & \textsc{957g43f} & \textsc{1000g} & \textsc{1000g} & \textsc{1000g} \\[0.25ex]
			\Qext{\sqrt{-3}} & & \textsc{904g96f} & \textsc{1000g} & \textsc{1000g} & \textsc{1000g} & & \textsc{924g76f} & \textsc{1000g} & \textsc{1000g} & \textsc{1000g} \\[0.5ex]
			\Qext{\sqrt{-7}} & & \textsc{939g61f} & \textsc{1000g} & \textsc{1000g} & \textsc{1000g} & & \textsc{961g39f} & \textsc{1000g} & \textsc{1000g} & \textsc{1000g} \\[0.5ex]
			\Qext{\sqrt{-11}} & & \textsc{979g21f} & \textsc{1000g} & \textsc{999g1f} & \textsc{1000g} & & \textsc{975g2b23f} & \textsc{1000g} & \textsc{995g5f} & \textsc{1000g} \\[0.5ex]

			\hline
		\end{tabular}
		\label{tbl:complex_quadratic-gamma_1-extensions}
	\end{table}

	Observe that when testing the field \( \Qext{-11} \) with complex \( C_i \) and \( \gamma = 2.0 \) the results were slightly worse than when \( \gamma = \gamma_1 \).
	This is likely because for this field \( \gamma_1 = \sqrt{22}/2 > 2 \), so \( \gamma = 2.0 \) is too small to satisfy the required constraints in Section~\ref{sub:algorithm_overview}.
	This supposition is strengthened by the observation that when \( \gamma = 3.0 > \sqrt{22}/2 \) the results are good again.
	
	We note a couple of \BAD{} results for the reduction method with complex \( C_i \). 
	In none of these cases did \APSLQ{} produce anything but a \GOOD{} result (if we ignore the case described in the previous paragraph).
	Nonetheless one or two bad results out of a pool of one thousand is hardly a poor result.
	
	The cases where \( D \in \{-5,-6,-10\} \) are cases where \( \gamma_1 \) does not exist and so the three conditions (\ref{cnd:rho}), (\ref{cnd:tau_rho}), and (\ref{cnd:tau_rho_gamma}) from Section~\ref{sub:algorithm_overview} are not satisfied.
	These results are summarised in Table~\ref{tbl:complex_quadratic-no_gamma_1-extensions}.
		
	\begin{table}
		\centering
		\caption{Complex quadratic fields without \( \gamma_1 \)}
		\begin{tabular}{ccccccccc}
			\hline
			\multirow{3}{*}{Field} & \hspace{2ex} & \multicolumn{3}{c}{Small Coefficients} & \hspace{2ex} & \multicolumn{3}{c}{Large Coefficients} \\
			& & \multirow{2}{*}{Reduction} & \multicolumn{2}{c}{Algebraic \PSLQ} & & \multirow{2}{*}{Reduction} & \multicolumn{2}{c}{Algebraic \PSLQ} \\
			& & & \( \gamma = 2.0 \) & \( \gamma = 3.0 \) & & & \( \gamma = 2.0 \) & \( \gamma = 3.0 \)\\
			\hline
			\\
			\multicolumn{9}{l}{Real \(C_i\)} \\
			\hline
			\rule{0pt}{2.5ex}\Qext{\sqrt{-5}} & & \textsc{994g6f} & \textsc{997g3f} & \textsc{1000g} & & \textsc{1000g} & \textsc{983g2b15f} & \textsc{992g2b6f} \\[0.25ex]
			\Qext{\sqrt{-6}} & & \textsc{996g4f} & \textsc{997g3f} & \textsc{999g1f} & & \textsc{998g2f} & \textsc{986g1b13f} & \textsc{996g1b3f} \\[0.5ex]
			\Qext{\sqrt{-10}} & & \textsc{1000g} & \textsc{999g1f} & \textsc{999g1f} & & \textsc{1000g} & \textsc{993g7f} & \textsc{994g6f} \\[0.5ex]
			\hline
			\\
			\multicolumn{9}{l}{Complex \(C_i\)} \\
			\hline
			\rule{0pt}{2.5ex}\Qext{\sqrt{-5}} & & \textsc{995g5f} & \textsc{158g842f} & \textsc{187g813f} & & \textsc{997g3f} & \textsc{164g836f} & \textsc{182g818f} \\[0.25ex]
			\Qext{\sqrt{-6}} & & \textsc{999g1f} & \textsc{136g864f} & \textsc{143g857f} & & \textsc{1000g} & \textsc{59g941f} & \textsc{60g940f} \\[0.5ex]
			\Qext{\sqrt{-10}} & & \textsc{1000g} & \textsc{40g960f} & \textsc{42g958f} & & \textsc{1000g} & \textsc{1000f} & \textsc{1000f} \\[0.5ex]
			\hline
		\end{tabular}
		\label{tbl:complex_quadratic-no_gamma_1-extensions}
	\end{table}	
	
	Observe that for real \( C_i \) the results are mostly good, despite the algorithm conditions not being satisfied.
	This is similar to the results for \( \Qext{\sqrt{-11}} \), highlighted above, that also failed those conditions.
	Contrast these to the cases with complex \( C_i \).
	
	The cases with complex \( C_i \) perform exceptionally poorly for \APSLQ{}. 
	This ought not be especially surprising since these fields do not satisfy the required conditions. 
	It is perhaps more remarkable that the results for the real \( C_i \) case are so good. 
	However, the reduction method gives consistently good results. 
	If we can find a way to reliably find correct algebraic integer relations from the incorrect ones often given by this method, it should prove to be remarkably robust. 
	% subsection complex_quadratic_extension_fields (end)
	
	% section experimental_results (end)

	\section{Further Work} % (fold)
	\label{sec:further_work}
		Further tests are being run which look more closely at the relationship between integer coefficient size, input vector size, and the precision necessary to find an integer relation.
		These tests also examine how the algorithm performs with problems consisting of extra constants than those that are known to be in the integer relation (i.e., relations with constants whose coefficient will be 0 in the integer relation).
		
		We suspect, based on some early proof-of-concept tests performed while implementing \APSLQ{}, that the reduction method will require more precision than \APSLQ{} for the same problem instance. 
		The above further tests should quantify that, if it is correct.
	
		Work is ongoing to find a theoretical framework with which to further modify the \APSLQ{} algorithm so that we may handle the real quadratic integer case, and the complex quadratic integer cases that do not satisfy the requirements from Section~\ref{sub:algorithm_overview}.
		
		Work is also ongoing to ascertain a method of reliably extracting algebraic integers in the complex quadratic reduction case. 
	% section further_work (end)

	% Print the bibliography:
	% 1st alternative is old way; the 2nd is biblatex (preferred).
	\ifx\printbibliography\undefined
		\bibliographystyle{spmpsci}
	    \bibliography{APSLQ}
	\else\printbibliography\fi
	
\end{document}